\documentclass[12pt]{article} 
\usepackage{amssymb,amsmath}
\usepackage{amsthm}
\usepackage{graphicx}
\newtheorem{theorem}{Theorem}[section]

\newtheorem{corollary}[theorem]{Corollary}

\theoremstyle{definition}

\theoremstyle{remark}

\begin{document}
\newcommand{\beq}{\begin{equation}} 
\newcommand{\eeq}{\end{equation}}
\newcommand{\zz}{\mathbb{Z}}
\newcommand{\pp}{\mathbb{P}} 
\newcommand{\nn}{\mathbb{N}}
\newcommand{\qq}{\mathbb{Q}}
\newcommand{\rr}{\mathbb{R}}
\newcommand{\cc}{\mathbb{C}}
\newcommand{\parn}{\mathrm{Par}(n)}
\newcommand{\bm}[1]{{\mbox{\boldmath $#1$}}}
\newcommand{\con}{\mathrm{Comp}(n)}
\newcommand{\sn}{\mathfrak{S}_n} 
\newcommand{\fs}{\mathfrak{S}}
\newcommand{\st}{\,:\,} 
\newcommand{\as}{\mathrm{as}}
\newcommand{\is}{\mathrm{is}}
\newcommand{\lgn}{\mathrm{len}}
\newcommand{\dis}{\displaystyle}
\newcommand{\bea}{\begin{eqnarray}}
\newcommand{\eea}{\end{eqnarray}}
\newcommand{\be}{\begin{enumerate}}
\newcommand{\ee}{\end{enumerate}}
\newcommand{\beas}{\begin{eqnarray*}}
\newcommand{\eeas}{\end{eqnarray*}}

\thispagestyle{empty}

\vskip 20pt
\begin{center}
{\large\bf Two Enumerative Results on Cycles of 
  Permutations}\footnote{This material is based upon 
  work supported by the National Science Foundation under Grant
  No.~0604423. Any opinions, findings and conclusions or
  recommendations expressed in this material are those of the author
  and do not necessarily reflect those of the National Science
  Foundation.} 
%  \footnote{2000 Mathematics
%  Subject Classification:  
%  05A15\\Key words and phrases: permutation, alternating sequence}
\vskip 15pt
{\bf Richard P. Stanley}\\[.1in]
{\it Department of Mathematics}\\ 
{\it Massachusetts Institute of Technology}\\
{\it Cambridge, MA 02139, USA}\\
{\texttt{rstan@math.mit.edu}}\\[.1in]
{\emph{In memory of Tom Brylawski}}\\[.2in]
{\bf\small version of 13 January 2009}\\
\end{center}

\begin{abstract}
Answering a question of B\'ona, it is shown that for $n\geq
2$ the probability that 1 and 2 are in the same cycle of a product of
two $n$-cycles on the set $\{1,2,\dots,n\}$ is $1/2$ if $n$ is odd and
$\frac 12-\frac{2}{(n-1)(n+2)}$ if $n$ is even. Another result
concerns the polynomial $P_\lambda(q)=\sum_w
q^{\kappa((1,2,\dots,n)\cdot w)}$, where $w$ ranges over all
permutations in 
the symmetric group $\sn$ of cycle type $\lambda$, $(1,2,\dots,n)$
denotes the $n$-cycle $1\rightarrow 2\rightarrow\cdots\rightarrow
n\rightarrow 1$, and $\kappa(v)$ denotes the number of cycles of the
permutation $v$. A formula is obtained for $P_\lambda(q)$ from which
it is deduced that all zeros of $P_\lambda(q)$ have real part 0.
\end{abstract}

\section{Introduction.} \label{sec1} 
Let $\lambda=(\lambda_1,\lambda_2,\dots)$ be a partition of $n$,
denoted $\lambda\vdash n$. In general, we use notation and terminology
involving partitions and symmetric functions from \cite[Ch.~7]{ec2}.
Let $\sn$ denote the symmetric group of all permutations of
$[n]=\{1,2,\dots,n\}$. If $w\in\sn$ then write $\rho(w)=\lambda$ if $w$
has cycle type $\lambda$, i.e., if the (nonzero) $\lambda_i$'s are the
lengths of the cycles of $w$. The conjugacy classes of $\sn$ are given
by $K_\lambda=\{ w\in\sn\st \rho(w)=\lambda\}$.

The ``class multiplication problem'' for $\sn$ may be stated as
follows. Given $\lambda,\mu,\nu\vdash n$, how many pairs $(u,v)\in
\sn\times \sn$ satisfy $u\in K_\lambda$, $v\in K_\mu$, $uv\in
K_\nu$? The case when one of the partitions is $(n)$ (i.e., one of
the classes consists of the $n$-cycles) is particularly interesting
and has received much attention. For a sample of some recent work, see
\cite{biane}\cite{irving} \cite{po-sc}.
% [goupil,schaeffer, jackson, boccara,  rs, ...] 
In this paper we make two contributions to this subject. For
the first, we solve a problem of B\'ona and Flynn \cite{bona} asking
what is the probability that two fixed elements of $[n]$ lie in the
same cycle of the product of two random $n$-cycles. In particular, we
prove the conjecture of B\'ona that this probability is $1/2$ when $n$
is odd.  Our method of proof is an ugly computation based on a formula
of Boccara \cite{boccara}. The technique can be generalized, and as an
example we compute the probability that three fixed elements of $[n]$
lie in the same cycle of the product of two random $n$-cycles.

For our second result, let $\kappa(w)$ denote the number of cycles of
$w\in \sn$, and let $(1,2,\dots,n)$ denote the $n$-cycle $1\rightarrow
2\rightarrow\cdots\rightarrow n\rightarrow 1$.  For $\lambda\vdash n$,
define the polynomial
  \beq P_\lambda(q)=\sum_{\rho(w)=\lambda} q^{\kappa((1,2,\dots,n)\cdot
    w)}. \label{eq:pdef} \eeq
In Theorem~\ref{thm3} we obtain a formula for $P_\lambda(q)$. We also
prove from this formula (Corollary~\ref{cor:iz}) that every zero of
$P_\lambda(q)$ has real part 0.

\section{A problem of B\'ona.} \label{sec2} Let $\pi_n$ denote the
probability that if two $n$-cycles $u,v$ are chosen uniformly at
random in $\sn$, then 1 and 2 (or any two elements $i$ and $j$ by
symmetry) appear in the same cycle of the product $uv$. Mikl\'os
B\'ona conjectured (private communication) that $\pi_n=1/2$ if $n$ is
odd, and asked about the value when $n$ is even. For the reason behind
this conjecture, see B\'ona and Flynn \cite{bona}. In this section we
solve this problem. Let us note that it is easy to see (a
straightforward generalization of \cite[Prop.~6.18]{bona2}) that the
probability that $1,2,\dots,k$ appear in the same cycle of
a random permutation in $\sn$ is $1/k$ for $k\geq n$.

\begin{theorem} \label{thm1}
For $n\geq 2$ we have
  $$ \pi_n = \left\{ \begin{array}{rl} \frac 12, & n\
      \mathrm{odd}\\[.5em] 
     \frac 12-\frac{2}{(n-1)(n+2)}, & n\ \mathrm{even}.
    \end{array} \right. $$
\end{theorem}

\proof
First note that if $w\in \sn$ has cycle type $\lambda$, then the
probability that 1 and 2 are in the same cycle of $w$ is
   $$ q_\lambda = \frac{\sum \binom{\lambda_i}{2}}{\binom n2} =
             \frac{\sum \lambda_i(\lambda_i-1)}{n(n-1)}. $$
Let $a_\lambda$ be the number of pairs $(u,v)$ of $n$-cycles in $\sn$
for which $uv$ has type $\lambda$. Then
   $$ \pi_n =\frac{1}{(n-1)!^2} \sum_{\lambda\vdash n} 
    a_\lambda q_\lambda. $$
By Boccara \cite{boccara} the number of ways to write a fixed
permutation $w\in \sn$ of type $\lambda$ as a product of two
$n$-cycles is 
   $$ (n-1)!\int_0^1\prod_i \left(x^{\lambda_i}-(x-1)^{\lambda_i}
        \right)dx. $$
Let $n!/z_\lambda$ denote the number of permutations $w\in \sn$ of
type $\lambda$. We get
  \beas \pi_n & = & \frac{1}{(n-1)!^2}\sum_{\lambda\vdash n}
    \frac{n!}{z_\lambda} 
   \left(\sum_i\frac{\lambda_i(\lambda_i-1)}{n(n-1)}\right)\\ &  &
    \qquad \cdot(n-1)!\int_0^1\prod_i
    \left(x^{\lambda_i}-(x-1)^{\lambda_i} \right)dx\\ & = &
   \frac{1}{n-1}\sum_{\lambda\vdash n}
    z_\lambda^{-1} 
   \left(\sum_i \lambda_i(\lambda_i-1)\right)
    \int_0^1\prod_i \left(x^{\lambda_i}-(x-1)^{\lambda_i}
        \right)dx. \eeas

Now let $p_\lambda(a,b)$ denote the power sum symmetric function
$p_\lambda$ in the two variables $a,b$, and let $\ell(\lambda)$ denote
the length (number of parts) of $\lambda$. It is easy to check that
  $$ 2^{-\ell(\lambda)}\left(\frac{\partial^2}{\partial a^2}-
       \frac{\partial^2}{\partial a\partial b}\right)
    \left. p_\lambda(a,b)\right|_{a=b=1} = \sum
    \lambda_i(\lambda_i-1). $$
By the exponential formula (permutation version)
\cite[Cor.~5.1.9]{ec2} or by \cite[Prop.~7.7.4]{ec2},
  $$ \sum_{n\geq 0} \sum_{\lambda\vdash n}
    z_\lambda^{-1} 2^{-\ell(\lambda)} 
   p_\lambda(a,b)
    \left(\prod_i \left(x^{\lambda_i}-(x-1)^{\lambda_i}
        \right)\right)t^n $$
  $$ \qquad =
   \exp \sum_{k\geq 1}\frac 1k \left(\frac{a^k+b^k}{2}\right)
         (x^k-(x-1)^k)t^k. $$
It follows that $(n-1)\pi_n$ is the coefficient of $t^n$ in
  $$ \hspace{-1in} F(t):= $$
  $$ \int_0^1 \left(\frac{\partial^2}{\partial a^2}-
       \frac{\partial^2}{\partial a\partial b}\right)
      \left. \exp\left[ \sum_{k\geq 1}\frac 1k
       \left(\frac{a^k+b^k}{2}\right) 
         (x^k-(x-1)^k)t^k\right]\right|_{a=b=1}\,dx. $$
We can easily perform this computation with Maple, giving
  \beas F(t) & = & \int_0^1 \frac{t^2(1-2x-2tx+2tx^2)}{(1-t(x-1))
     (1-tx)^3}dx\\
    & = & \frac{1}{t^2}\log(1-t^2)+\frac 32+\frac{-\frac
    12+t}{(1-t)^2}. \eeas
Extract the coefficient of $t^n$ and divide by $n-1$ to
obtain $\pi_n$ as claimed. 
\qed

\medskip
It is clear that the argument used to prove Theorem~2 can be
generalized. For instance, using the fact that
  $$  3^{-\ell(\lambda)-1}\left(\frac{\partial^3}{\partial a^3}-
       3\frac{\partial^3}{\partial a^2\partial b}
     +2\frac{\partial^3}{\partial a\partial b\partial c}\right)
    \left. p_\lambda(a,b,c)\right|_{a=b=c=1} $$
   $$ \qquad = \sum
    \lambda_i(\lambda_i-1)(\lambda_i-2), $$
we can obtain the following result. 

\begin{theorem} \label{thm2}
Let $\pi_n^{(3)}$ denote the probability
  that if two $n$-cycles $u,v$ are chosen uniformly at random in
  $\sn$, then $1, 2$, and $3$ appear in the same cycle of the product
  $uv$. Then for $n\geq 3$ we have
  $$ \pi_n^{(3)} = \left\{ \begin{array}{rl}
     \frac 13+\frac{1}{(n-2)(n+3)}, & n\ \mathrm{odd}\\[.5em]
    \frac 13 -\frac{3}{(n-1)(n+2)}, & n\ \mathrm{even}. \end{array}
     \right. $$ 
\end{theorem}

Are there simpler proofs of Theorems~\ref{thm1} and \ref{thm2},
especially Theorem~\ref{thm1} when $n$ is odd?

\section{A polynomial with purely imaginary zeros} \label{sec3}

Given $\lambda\vdash n$, let $P_\lambda(q)$ be defined by
equation~\eqref{eq:pdef}.  Let $(a)_n$ denote the falling factorial
$a(a-1)\cdots (a-n+1)$. Let $E$ be the backward shift operator on
polynomials in $q$, i.e., $Ef(q)=f(q-1)$. 

\begin{theorem} \label{thm3}
Suppose that $\lambda$ has length $\ell$. Define the polynomial  
   $$ g_\lambda(t)=\frac{1}{1-t}\prod_{j=1}^\ell
     (1-t^{\lambda_j}). $$ 
Then 
  \beq P_\lambda(q) = z_\lambda^{-1} g_\lambda(E)(q+n-1)_n. 
    \label{eq:plaq} \eeq
\end{theorem}

\proof
Let $x=(x_1,x_2,\dots)$, $y=(y_1,y_2,\dots)$, and $z=(z_1,z_2,\dots)$
be three disjoint sets of variables. Let $H_\mu$ denote the
product of the hook lengths of the partition $\mu$ (defined e.g.\
in \cite[p.~373]{ec2}). Write $s_\lambda$ and $p_\lambda$ for the Schur
function and power sum symmetric function indexed by $\lambda$. The
following identity is the case $k=3$ of 
\cite[Prop.~2.2]{h-s-s} and \cite[Exer.~7.70]{ec2}:
  \beq \sum_{\mu\vdash n}H_\mu s_\mu(x)s_\mu(y)s_\mu(z) =
    \frac{1}{n!}\sum_{uvw=1\,\mathrm{in}\,\sn} p_{\rho(u)}(x)
    p_{\rho(v)}(y)p_{\rho(w)}(z). \label{eq:fund} \eeq
For a symmetric function $f(x)$ let $f(1^q)=f(1,1,\dots,1,0,0,\dots)$
($q$ 1's). Thus $p_{\rho(w)}(1^q)=q^{\kappa(w)}$. Let
$\chi^\lambda(\mu)$ denote the irreducible character of 
$\sn$ indexed by $\lambda$ evaluated at a permutation of cycle type
$\mu$ \cite[{\S}7.18]{ec2}. Recall \cite[Cor.~7.17.5 and
Thm.~7.18.5]{ec2} that 
   $$ s_\mu=\sum_{\nu\vdash n}z_\nu^{-1}\chi^\mu(\nu)p_\nu, $$
where $\#K_\nu=n!/z_\nu$ as above. Take the coefficient of 
$p_n(x)p_\lambda(y)$ in equation~\eqref{eq:fund} and set
$z=1^q$. Since there are $(n-1)!$ $n$-cycles $u$, the right-hand side
becomes $\frac 1n P_\lambda(q)$. Hence
  \beq P_\lambda(q) = n\sum_{\mu\vdash n}H_\mu z_n^{-1}\chi^\mu(n)
     z_\lambda^{-1}\chi^\mu(\lambda)s_\mu(1^q). \label{eq:int} \eeq
Write $\sigma(i)=\langle n-i,1^i\rangle$, the ``hook'' with one part
equal to $n-i$ and $i$ parts equal to $1$, for $0\leq i\leq n-1$. Now
$z_n=n$, and e.g.\ by \cite[Exer.~7.67(a)]{ec2} we have
   $$ \chi^\mu(n) = \left\{ \begin{array}{rl} 
    (-1)^i, & \mathrm{if}\ \mu=\sigma(i),\ 0\leq i\leq n-1\\[.1in]
    0, & \mathrm{otherwise}. \end{array} \right. $$
Moreover, $s_{\sigma(i)}(1^q)=(q+n-i-1)_n H_{\sigma(i)}^{-1}$
by the hook-content formula \cite[Cor.~7.21.4]{ec2}. Therefore we get
from equation~\eqref{eq:int} that
  \beq P_\lambda(q) = z_\lambda^{-1}\sum_{i=0}^{n-1}(-1)^i
  \chi^{\sigma(i)}(\lambda)(q+n-i-1)_n. \label{eq:int2} \eeq
\indent The following identity is a simple consequence of Pieri's rule
\cite[Thm.~7.15.7]{ec2} and appears in \cite[I.3, Ex.~14]{macd}:
  $$ \prod_i \frac{1+tx_i}{1-ux_i} = 1+(t+u)\sum_{i=0}^{n-1} 
    s_{\sigma(i)}t^iu^{n-i-1}. $$
Substitute $-t$ for $t$, set $u=1$ and take the scalar product with
$p_\lambda$. Since $\langle s_\mu,p_\lambda\rangle =\chi^\mu(\lambda)$
the right-hand side becomes $(1-t)\sum_{i=0}^{n-1} 
(-1)^i\chi^{\sigma(i)}(\lambda)t^i$. On the other hand, the
left-hand side is given by 
  \beas \left\langle \exp\left( \sum_{n\geq 1}\frac{p_n}{n}\right)
   \cdot\exp\left(-\sum_{n\geq 1} \frac{p_n}{n}t^n\right),
    p_\lambda\right\rangle
    & = & \left\langle \exp\left( \sum_{n\geq 1}\frac{p_n}{n}(1-t^n)
    \right), p_\lambda\right\rangle\\ & = & 
   \prod_{i=1}^\ell \left(1-t^{\lambda_i}\right), \eeas
by standard properties of power sum symmetric functions
\cite[{\S}7.7]{ec2}. Hence
  $$ \sum_{i=0}^{n-1} (-1)^i\chi^{\sigma(i)}(\lambda)t^i =
           g_\lambda(t). $$ 
Comparing with equation~\eqref{eq:int2} completes the proof.
\qed 

\medskip
\textsc{Note.} Since $(1-E)(q+n)_{n+1}=(n+1)(q+n-1)_n$,
equation~\eqref{eq:plaq} can be rewritten as
  \beq P_\lambda(q) = \frac{1}{(n+1)z_\lambda}g'_\lambda(E)
      (q+n)_{n+1}, \label{eq:gprime} \eeq
where $g'_\lambda(t)=\prod_{j=1}^\ell (1-t^{\lambda_j})$.
\medskip
 
The zeros of the polynomial $P_\lambda(q)$ have an interesting
property that will follow from the following result.

\begin{theorem} \label{thm:iz}
Let $g(t)$ be a complex polynomial of degree exactly $d$, such that
every zero of $g(t)$ lies on the circle $|z|=1$. Suppose that the
multiplicity of 1 as a root of $g(t)$ is $m\geq 0$. Let 
$P(q)=g(E)(q+n-1)_n$.
 \be\item[(a)] If $d\leq n-1$, then 
   $$ P(q) = (q+n-d-1)_{n-d}\,Q(q), $$
where $Q(q)$ is a polynomial of degree $d-m$ for which every zero has
real part $(d-n+1)/2$. 
   \item[(b)] If $d\geq n-1$, then $P(q)$ is a polynomial of degree
$n-m$ for which every zero has real part $(d-n+1)/2$. 
  \ee
\end{theorem}

\proof
First, the statements about the degrees of $Q(q)$ and $P(q)$ are
clear; for we can write $g(t)=c\prod_u(t-u)$ and apply the factors
$t-u$ consecutively. If $h(q)$ is any polynomial and $u\neq 1$ then
$\deg\, (E-u)h(q)=\deg h(q)$, while $\deg\, (E-1)h(q)= \deg h(q) -1$.

The remainder of the proof is by induction on $d$. The base case $d=0$
is clear. Assume the statement for $d<n-1$. Thus for $\deg g(t) =d$ we
have
  \beas g(E)(q+n-1)_n & = & (q+n-d-1)_{n-d}\,Q(q)\\ & = & 
     (q+n-d-1)_{n-d}\prod_j\left(q-\frac{d-n+1}{2}
     -\delta_j i\right)  \eeas
for certain real numbers $\delta_j$. Now
   $$ \hspace{-1in}(E-u)g(E)(q+n-1)_n $$
   \vspace{-3em}
   \beas & = &
    (q+n-d-1)_{n-d}\,Q(q)-u(q+n-d-2)_{n-d}\,Q(q-1)\\ & = &
     (q+n-d-2)_{n-d-1}[(q+n-d-1)Q(q)-u(q-1)\,Q(q-1)]\\ & = &
      (q+n-d-2)_{n-d-1}Q'(q), \eeas
say. The proof now follows from a standard argument (e.g.,
\cite[Lemma~9.13]{po-st}), which we give for the sake of
completeness. Let 
$Q'(\alpha+\beta i)=0$, where $\alpha,\beta\in\rr$.  Thus
  $$ (\alpha+\beta i +n-d-1)\prod_j
     \left(\alpha+\beta i-\frac{d-n+1}{2}-\delta_j i\right) $$ 
  $$ = u(\alpha+\beta i -1)\prod_j
     \left(\alpha-1+\beta i-\frac{d-n+1}{2}-\delta_j i\right). $$
Letting $|u|=1$ and taking the square modulus gives
  $$ \frac{(\alpha+n-d-1)^2+\beta^2}{(\alpha-1)^2+\beta^2}
    \prod_j\frac{\left(\alpha-
    \frac{d-n+1}{2}\right)^2+(\beta-\delta_j)^2}
     {\left(\alpha-1-
    \frac{d-n+1}{2}\right)^2+(\beta-\delta_j)^2} = 1. $$
If $\alpha<(d-n+2)/2$ then 
   $$ (\alpha+n-d-1)^2-(\alpha-1)^2<0 $$
and
  $$ \left( \alpha-\frac{d-n+1}{2}\right)^2
     < \left( \alpha-1-\frac{d-n+1}{2}\right)^2. $$
The inequalities are reversed if $\alpha>(d-n+2)/2$. Hence $\alpha = 
(d-n+2)/2$, so the theorem is true for $d\leq n-1$. 

For $d\geq n-1$ we continue the induction, the base case now being
$d=n-1$ which was proved above. The induction step is completely
analogous to the case $d\leq n-1$ above, so the proof is complete.
\qed

\begin{corollary} \label{cor:iz} 
The polynomial $P_\lambda(q)$ has degree $n-\ell(\lambda)+1$, and
every zero of $P_\lambda(q)$ has real part 0. 
\end{corollary}

\proof
The proof is immediate from Theorem~\ref{thm3} and the special case
$g(t)=g_\lambda(t)$ (as defined in Theorem~\ref{thm3}) and $d=n-1$ of
Theorem~\ref{thm:iz}.
\qed
\medskip

It is easy to see from Corollary~\ref{cor:iz} (or from considerations
of parity) that $P_\lambda(q)=(-1)^n P_\lambda(-q)$. Thus we can write
  $$ P_\lambda(q)=\left\{ \begin{array}{rl} R_\lambda(q^2), & n\
      \mathrm{even}\\[.1in]  qR_\lambda(q^2), & n\ \mathrm{odd}, 
    \end{array} \right. $$
for some polynomial $R_\lambda(q)$. It follows from
Corollary~\ref{cor:iz} that $R_\lambda(q)$ has (nonpositive) real
zeros. In particular (e.g., \cite[Thm.~2]{unim}) the coefficients of
$R_\lambda(q)$ are log-concave with no external zeros, and hence
unimodal.

The case $\lambda=(n)$ is especially interesting. Write $P_n(q)$ for
$P_{(n)}(q)$. From equation~\eqref{eq:gprime} we have
 $$ P_n(q) = \frac{1}{n(n+1)}((q+n)_{n+1}-(q)_{n+1}). $$
Now
   $$ (q)_{n+1} = (-1)^{n+1} (-q+n)_{n+1} $$ 
and
   $$ (q+n)_{n+1} = \sum_{k=1}^{n+1} c(n+1,k)q^k, $$
where $c(n+1,k)$ is the signless Stirling number of the first kind
(the number of permutations $w\in\fs_{n+1}$ with $k$ cycles)
\cite[Prop.~1.3.4]{ec1}. Hence 
  $$ \frac{1}{n(n+1)}((q+n)_{n+1}-(q)_{n+1}) =
   \frac{1}{\binom{n+1}{2}} 
         \sum_{k\equiv n\,(\mathrm{mod}\,2)} c(n+1,k)x^k. $$
We therefore obtain the following result.

\begin{corollary} \label{cor:cnk}
The number of $n$-cycles $w\in\sn$ for which $w\cdot(1,2,\dots,n)$ has
exactly $k$ cycles is 0 if $n-k$ is odd, and is otherwise equal to
$c(n+1,k)/\binom{n+1}{2}$. 
\end{corollary}

Is there a simple bijective proof of Corollary~\ref{cor:cnk}?

Let $\lambda,\mu\vdash n$. A natural generalization of $P_\lambda(q)$
is the polynomial 
   $$ P_{\lambda,\mu}(q)=\sum_{\rho(w)=\lambda} q^{\kappa(w_\mu\cdot
    w)}, $$
where $w_\mu$ is a fixed permutation in the conjugacy class $K_\mu$.
Let us point out that it is \emph{false} in general that every zero
of $P_{\lambda,\mu}(q)$ has real part 0. For instance, 
  $$ P_{332,332}(q) = q^8+35q^6+424q^4+660q^2, $$
four of whose zeros are approximately $\pm 1.11366\pm 4.22292i$.

%\pagebreak

\end{document}